\documentclass[pdflatex,sn-mathphys-num]{sn-jnl}
\usepackage{graphicx}%
\usepackage{multirow}%
\usepackage{amsmath,amssymb,amsfonts}%
\usepackage{amsthm}%
\usepackage{mathrsfs}%
\usepackage[title]{appendix}%
\usepackage{xcolor}%
\usepackage{textcomp}%
\usepackage{manyfoot}%
\usepackage{booktabs}%
\usepackage{algorithm}%
\usepackage{algorithmicx}%
\usepackage{algpseudocode}%
\usepackage{listings}%
\usepackage[T1]{fontenc}
\usepackage[utf8]{inputenc}
\usepackage{fontawesome5}
\newcommand{\jumplink}[1]{%
\href{https://wuprover.github.io/groebner_proj/docs/find/?pattern=#1\#src}{\small\faExternalLink*}%
}
\usepackage{tikz}
\usetikzlibrary{shapes, arrows.meta, positioning} 
\usepackage{wrapfig} 
\usepackage{lipsum}
\usetikzlibrary{arrows.meta, positioning, shapes.geometric, calc}
\tikzset{
    def_node/.style={
        rectangle, 
        rounded corners=3pt, 
        draw=teal!80, 
        fill=teal!10, 
        very thick, 
        minimum width=3cm, 
        minimum height=0.8cm, 
        align=center, 
        font=\sffamily\footnotesize
    },
    thm_node/.style={
        rectangle, 
        draw=blue!80, 
        fill=blue!10, 
        very thick, 
        minimum width=3cm, 
        minimum height=0.8cm, 
        align=center, 
        font=\sffamily\footnotesize
    },
    my_arrow/.style={-{Stealth[scale=1.0]}, thick, draw=gray!70}
}
\definecolor{keywordcolor}{rgb}{0.7, 0.1, 0.1}   
\definecolor{symbolcolor}{rgb}{0.0, 0.1, 0.6}    
\definecolor{sortcolor}{rgb}{0.1, 0.5, 0.1}      

\lstnewenvironment{leancode}{\lstset{frame=single,language=lean,abovecaptionskip=-\medskipamount,basicstyle={\ttfamily\footnotesize}}}{}
\newcommand{\lean}[1]{\texttt{#1}}
\newtheorem{theorem}{Theorem}
\newtheorem{lemma}{Lemma}
\newtheorem{proposition}[theorem]{Proposition}%
\newtheorem*{runningexample}{Example 1}
\newtheorem{definition}{Definition}%
\usepackage{diagbox}
\usepackage{xifthen}
\newcommand*{\mathlibdef}[2][]{\href{https://leanprover-community.github.io/mathlib4_docs/find/?pattern=#2\#src}{\ifthenelse{\isempty{#1}}{\lean{#2}}{#1}}}
\newcommand*{\ourdef}[2][]{\href{https://wuprover.github.io/groebner_proj/docs/find/?pattern=#2\#src}{\ifthenelse{\isempty{#1}}{\lean{#2}}{#1}}}
\newcommand*{\leandef}[2][]{\mathlibdef[#1]{#2}}
\newcommand{\M}{\ensuremath{\mathcal{M}}}
\DeclareMathOperator\LM{LM}
\DeclareMathOperator\LT{LT}
\DeclareMathOperator\LC{LC}
\DeclareMathOperator\multideg{deg}
\DeclareMathOperator\degree{\multideg}

\DeclareMathOperator\spoly{S}

\newcommand\kxi{k[x_i]_{i\in\sigma}}
\newcommand\kxiv{k[x_i]_{i\in\sigma'}}
\newcommand\leantodo[1]{}
\usepackage{caption}
\begin{document}
\title[Formalizing Gröbner Basis Theory in Lean ]{Formalizing Gröbner Basis Theory in Lean}
\author[1]{\fnm{Junyu } \sur{Guo}}\email{guojy228@mail2.sysu.edu.cn}
\equalcont{Equal contribution. We list the first authors alphabetically by last name.}
\author[2]{\fnm{Hao} \sur{Shen}}\email{shenhao24@amss.ac.cn}
\equalcont{Equal contribution. We list the first authors alphabetically by last name.}
\author[2]{\fnm{Junqi} \sur{Liu}}\email{liujunqi@amss.ac.cn}
\author[2]{\fnm{Lihong} \sur{Zhi}}\email{lzhi@mmrc.iss.ac.cn}
\affil[1]{\orgdiv{Institute of Logic and Cognition, Department of Philosophy}, \orgname{Sun Yat-sen University}, \orgaddress{\city{Guangzhou}, \postcode{510275}, \state{Guangdong}, \country{China}}}
\affil[2]{\orgdiv{State Key Laboratory of Mathematical Sciences}, \orgname{Academy of Mathematics and Systems Science, University of Chinese Academy of
Science}, \orgaddress{\street{No.55
Zhongguancun East Road}, \city{Beijing}, \postcode{100190}, \state{Beijing}, \country{China}}}
\abstract{
We present a formalization of Gröbner basis theory in Lean~4, built on top of Mathlib's infrastructure for multivariate polynomials and monomial orders. Our development covers the core foundations of  Gröbner basis theory, including polynomial division with remainder, Buchberger’s criterion, and the existence and uniqueness of reduced Gröbner bases. We develop the theory uniformly for polynomial rings indexed by arbitrary types, enabling the treatment of Gröbner bases in rings with infinitely many variables. Furthermore, we connect the finite and infinite settings by showing that infinite-variable reduced Gröbner bases can be characterized via reduced Gröbner bases on finite-variable subrings through monomial-order embeddings and filter-based limit constructions.}
\maketitle




\section{Introduction}
Gröbner bases, introduced by Buchberger~\cite{buchberger1965algorithmus}, are a cornerstone of computational commutative algebra and algebraic geometry. They provide a unified algorithmic framework for fundamental problems such as ideal membership testing, elimination theory, and solving systems of polynomial equations. As a result, they serve as a central engine for algebraic reasoning, supporting automated deduction in polynomial algebra and applications such as cryptography, robotics, and control theory. Consequently, a comprehensive development of Gröbner basis theory is essential to advance formalized mathematics and certify the correctness of algebraic algorithms.
The most extensive formalizations of Gröbner basis theory to date have been developed in the Isabelle/HOL proof assistant~\cite{nipkow2002isabelle},  including certified implementations of advanced methods such as the Faug\`{e}re’s $F_4$ algorithm~\cite{faugere1999new} and signature based techniques~\cite{maletzky2021generic, maletzky2024grobner,maletzky2018grobner}.  
In Coq~\cite{bertot2013interactive}, early efforts focused on formalizing Buchberger's algorithm and establishing the constructive existence of Gröbner bases~\cite{ persson2001integrated,Thery2001}.
Additional efforts include formalizations  in ACL2~\cite{medina2010verified} and Mizar~\cite{schwarzweller2005grobner}.
Dermitzakis implemented S-polynomials and basic Gr\"{o}bner basis constructions in Lean~3~\cite{dermitzakis2019implementation}. While verifying the feasibility of such a development, this work remained largely isolated from the core Mathlib algebraic infrastructure. Poulsen and Guo explored an initial formalization of the Gröbner basis theory, respectively, within Lean~3~\cite{andreas2023gb} and the new Lean~4 ecosystem~\cite{HagbLeanGroebner}. However, they neither continue to maintain their formalization nor submit their work to Mathlib.
The Lean tactic \texttt{polyrith} has successfully used Gröbner basis techniques for proof automation. The tactic \texttt{polyrith} relies on SageMath~\cite{sagemath} to compute certificates, which are then verified inside Lean.  The Lean~4 tactic \texttt{grobner} provides Gröbner-basis-based solving for polynomial goals and is implemented as a thin wrapper around the general automation tactic \texttt{grind}. Although effective for goal-directed automation, they assist in verification rather than in providing a formal development of the Gröbner basis theory.
In contrast, our work develops the core foundational infrastructure for the Gröbner basis theory in Lean~4. Parts of this development have already been merged into Mathlib, and the formalization is designed to support long-term reuse at the library-level. A key distinguishing feature of our work is its level of generality. Unlike classical formalizations restricted to finite Gr\"{o}bner basis, we reuse Mathlib's general infrastructure, handle polynomial rings indexed by elements in any given types, and also support infinite Gröbner bases in the sense of Iima and Yoshino~\cite{IimaYoshino08}.
Beyond definitional generality, we formalize essential results bridging Gröbner bases over different polynomial rings, including a characterization of infinite reduced Gr\"{o}bner bases via filter-based limit of finite reduced Gr\"{o}bner bases. This bridges the finite and infinite settings within a single formal framework, providing a reusable infrastructure for future formalization of computational algebraic geometry.
\subsection*{Organization}
The remainder of this paper is organized as follows. Section~\ref{sec:preliminaries} reviews the mathematical preliminaries and existing Mathlib infrastructure, including multivariate polynomials and monomial orders. In
Section~\ref{sec:formalization}, we formalize the core parts of Gr\"{o}bner basis theory. We begin by introducing a degree with bottom element to resolve arithmetic ambiguities regarding the zero polynomial, followed by the formalization of polynomial division and the existence of remainders. We then define the Gröbner bases and reduced Gröbner bases, establish their fundamental properties, and formally verify Buchberger’s criterion.
Section~\ref{infinite_to_finite} introduces embeddings of monomial orders, proves invariance of Gröbner bases under these embeddings, and establishes the characterization of infinite Gröbner bases via limit constructions. Finally, Section~\ref{future} discusses future work on the certification of external computations.
\section{Preliminaries}\label{sec:preliminaries}
Lean is an open source interactive theorem prover based on dependent type theory with a small trusted kernel~\cite{de2015lean}. It integrates an interactive proof language with an efficient compiled programming language and a powerful metaprogramming framework. Its latest version is Lean~4~\cite{moura2021lean}. Within the Lean ecosystem, Mathlib provides a large and coherent library of formalized mathematics~\cite{mathlib}.
Mathlib provides foundations for polynomial algebra, including \mathlibdef{MvPolynomial}, \mathlibdef{Finsupp}, and recent developments formalizing monomial orders and multivariate polynomial division. Our formalization builds directly on this algebraic hierarchy, reusing core abstractions while remaining definitionally compatible with the library, thereby providing the infrastructure required for Gröbner basis theory.
Let $\sigma$ be a type of indices, and let 
\begin{equation}
 \label{Mset}
   \M:= (\sigma \to_0 \mathbb{N})    
\end{equation}
denote the type of finitely supported functions from $\sigma$ to $\mathbb{N}$. A monomial order on $\M$ is defined as follows.
 
\begin{definition}[{\leandef[Monomial Order]{MonomialOrder}}]
\label{monomial_order}
A monomial order is a binary relation $\le$ on $\M$ satisfying the following conditions:
\begin{enumerate}
    \item[(1)] 
    The relation $\le$ is a well-order on $\M$, i.e., it is a linear order and admits no infinite strictly descending chain   $
    a_1 > a_2 > a_3 > \cdots. $
    
    \item[(2)] For all $a, b, c \in \M$, if $a \le b$, then $a + c \le b + c$.
    
    \item[(3)] The zero element $0$ is the smallest element, i.e., $0 \le a$ for all $a \in \M$.
\end{enumerate}
\end{definition}
A standard example of such an ordering is lexicographic order. Assume that the index type $\sigma$ is equipped with a linear order. The corresponding lexicographic order $\le_{\operatorname{lex}}$ on $\M$ is defined as follows:
 
\begin{definition}[{\leandef[Lexicographic Order]{MonomialOrder.lex}}]
\label{lex_order}
Assume $\sigma$ is equipped with a linear order
 $\geq$, where $>$ is well-founded, i.e., it admits no infinite strictly ascending sequence $a_1 < a_2 < a_3 < \cdots $. 
For distinct $a, b \in \M$, we define $a <_{\operatorname{lex}} b$ if
\[
a(i) < b(i), \quad \text{where } i = \min \{ k \in \sigma \mid a(k) \neq b(k) \}.
\]
Since $a$ and $b$ have finite support, this minimum exists. 
\end{definition}
Let $f = \sum_{j \in \M} a_{j} x^{j}$ be a polynomial in $\kxi$, where only finitely many coefficients $a_j$ are nonzero, $j$ is called an exponent vector and for any $j$ such that $a_{j} \neq 0$, the expression $a_{j} x^{j}$ is called a term of $f$. 
 
Once a monomial order is fixed, we can define the leading term, leading coefficient and the leading monomial of $f$. These notions play a fundamental role in the Gröbner basis theory.
\begin{definition}[{\ourdef[Leading Term]{MonomialOrder.leadingTerm}}]
\label{leadingterm}
For a nonzero multivariate polynomial $f$, its leading term, denoted by $\LT(f)$, is the term of $f$ whose monomial is maximal with respect to the fixed monomial order. For the zero polynomial, the leading term is defined to be $0$.
\end{definition}
\begin{definition}[{\leandef[Leading Coefficient]{MonomialOrder.leadingCoeff}}]
\label{leadingcoeff}
The leading coefficient of a polynomial $f$, denoted by $\LC(f)$, is defined as the coefficient of the leading term of $f$.
\end{definition}
\begin{definition}[Leading Monomial]
Let a monomial order be fixed, the leading monomial of a polynomial $f$, denoted by $\LM(f)$, is defined as the monomial corresponding to the maximal exponent vector in the support of $f$ with respect to the monomial order.  If $f = 0$, then $\LM(f) = 0$.
\end{definition}
Although the leading monomial is not provided as a dedicated definition in Lean, it can be represented using the expression \texttt{monomial (m.degree f) 1} if $f \ne 0$. We postpone the formalization of the leading term to the next section. 
\section{Formalizing Gröbner Basis Theory in Lean}\label{sec:formalization}
\label{Formalizing_gb}
\begin{figure}[t]
    \centering
\begin{tikzpicture}[
    scale=0.85, transform shape,
    node distance=1.3cm and 1.3cm,
    definition/.style={
        rectangle, 
        rounded corners=3pt, 
        draw=teal!80, 
        fill=teal!10, 
        very thick, 
        minimum width=3.5cm, 
        minimum height=1cm, 
        align=center, 
        font=\sffamily\small
    },
    theorem/.style={
        rectangle, 
        draw=blue!80, 
        fill=blue!10, 
        very thick, 
        minimum width=3.5cm, 
        minimum height=1cm, 
        align=center, 
        font=\sffamily\small
    },
    arrow/.style={-{Stealth[scale=1.2]}, thick, draw=gray!80},
    label_text/.style={font=\scriptsize\itshape, color=gray}
]
\node[definition] (def1) {Definition \ref{monomial_order} \\Monomial Order};
\node[definition] (def3) [below left=1cm and 0.5cm of def1] {Remainder with \\ Bottom Element};
\node[definition] (def6) [below right=1cm and 0.5cm of def1] {Definition \ref{leadingterm}\\Leading Term};
\node[definition] (def5) [below=2.5cm of def1] {Definition \ref{isremainder}\\Remainder};
\node[theorem] (th1) [right=of def5] {Theorem \ref{remainder_exist}\\Remainder Existence};
\node[definition] (def7) [below=1.5cm of def5] {Definition \ref{gb}\\Gröbner Basis(GB)};
\node[theorem] (th2) [left=2cm of def7] {Theorem \ref{exist_gb}\\Existence of GB};
\node[theorem] (th4) [below=0.5cm of th2] {Proposition \ref{ideal_eq_span}\\$I = \langle G \rangle$};
\node[theorem] (th5) [below=0.5cm of th4] {Proposition \ref{gb_characterization_via_degree}\\GB Characterization\\via Degree};
\node[theorem] (th3) [right=2cm of def7] {Theorem \ref{ideal_mem}\\$p \in I \iff p \xrightarrow{G} 0$};
\node[theorem] (th6) [below=0.5cm of th3] {Proposition \ref{gb_remainder}\\GB Characterization\\via Remainder};
\node[theorem] (th8) [below=0.5cm of th6] {Theorem \ref{unique_of_remainder}\\Uniqueness of\\Remainder};
\node[definition] (def8) [below=4.5cm of def7] {Definition \ref{spoly}\\S-Polynomial};
\node[theorem] (th9) [right=of def8] {Theorem \ref{buchberger}\\Buchberger's Criterion};
\node[definition] (def10) [left=of def8] {Definition \ref{reduced_gb}\\Reduced GB};
\node[theorem] (th10) [below=of def10] {Theorem \ref{uniqueExists_reduced_gb}\\Unique Reduced GB};
\draw[arrow] (def1) -- (def3);
\draw[arrow] (def1) -- (def6);
\draw[arrow] (def3) -- (def5);
\draw[arrow] (def6) -- (def5); 
\draw[arrow] (def5) -- (th1);
\draw[arrow] (def5) -- (def7);
\draw[arrow] (def6) -- (def7);
\draw[arrow] (def7) -- (th2);
\draw[arrow] (def7) -- (th4);
\draw[arrow] (def7) -- (th5);
\draw[arrow] (def7) -- (th3);
\draw[arrow] (def7) -- (th8);
\draw[arrow, dashed] (th3) -- (th6); 
\draw[arrow] (def7) -- (th9); 
\draw[arrow] (def8) -- (th9); 
\draw[arrow] (def7) -- (def10); 
\draw[arrow] (def10) -- (th10); 
\draw[arrow, dashed] (th1) -- (th8); 
\end{tikzpicture}
\begin{centering}
\caption{Structure of the formalization of Gröbner basis}
\end{centering}
\label{sec3_overview}
\end{figure}
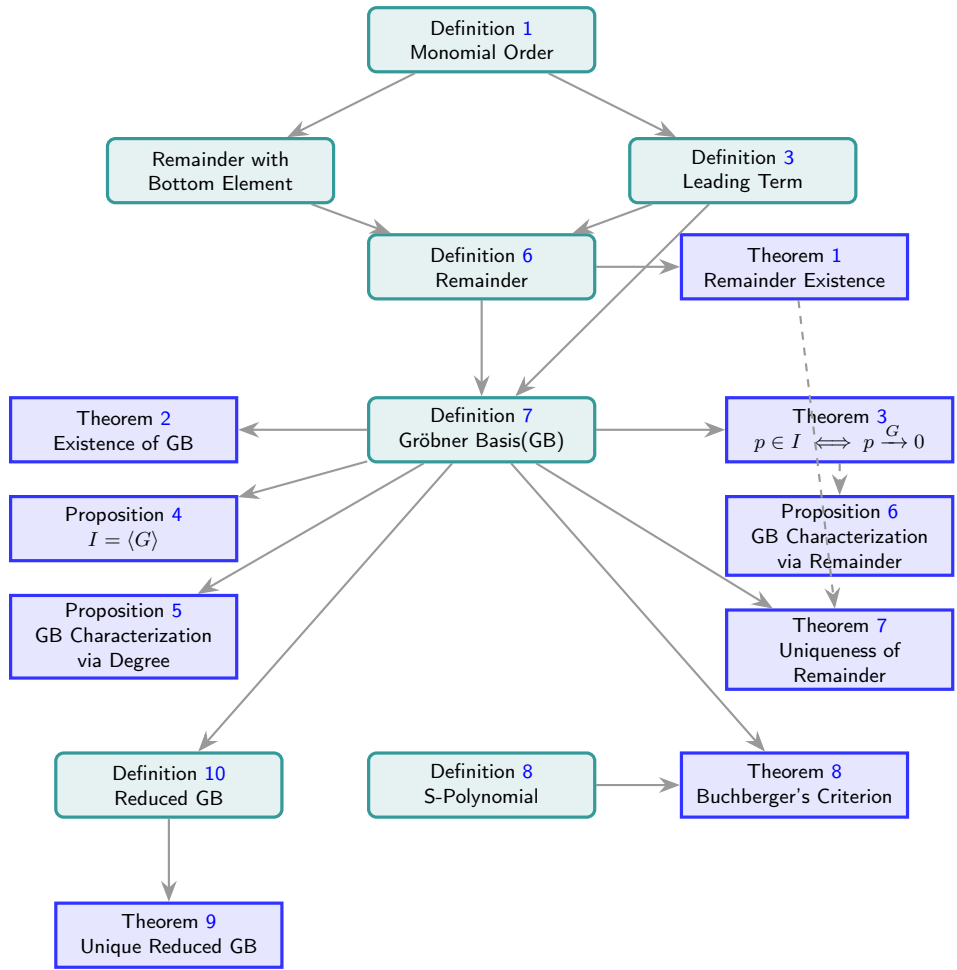

In this section, we present our formalization of Gröbner basis theory. While classical treatments typically consider polynomial rings in finitely many variables, the multivariate polynomial library in Mathlib is defined over an arbitrary index type $\sigma$, which is not necessarily finite. Accordingly, we develop Gröbner basis theory for polynomial rings with arbitrarily many variables. In the infinite-variable setting, the polynomial ring is generally not Noetherian, and ideals need not admit finite generating sets. To accommodate this, we extend the notion of Gröbner basis to allow infinite sets, and show that most fundamental properties continue to hold in this setting. Although the main theorems in this section are stated over fields, several results have also been formalized over more general commutative rings. An overview of the dependency structure of our formalization is given in Figure~\ref{sec3_overview}.
\subsection{Degree with Bottom Element}\label{subsec:degree}
\label{degree_with_bot}
Let $f = \sum_{j \in \M} a_j x^j$ be a non-zero polynomial. The degree of $f$, denoted by $\degree(f)$, is defined as
\[ \degree(f) := \max \{ j \in \M \mid a_j \neq 0 \}\]
where the maximum is taken with respect to the monomial order.  
In Mathlib, the monomial exponent type for multivariate polynomials is implemented as \texttt{Finsupp $\sigma$ $\mathbb{N}$}, written mathematically as $\sigma \to_0 \mathbb{N}$. This corresponds to the exponent type $\mathcal{M}$ introduced in~(\ref{Mset}).
According to the existing Mathlib convention, the degree of the zero polynomial is defined to be $0$. While convenient for computations on exponent vectors, this convention fails to preserve certain  properties of the degree function. 
In particular, the additive property 
\[
\deg(fg) = \deg(f) + \deg(g)
\]
does not hold without additional assumptions. For example, if $f = 0$ and $g$ have a positive degree, then $\deg(fg)=0$ while $\deg(f)+\deg(g)=\deg(g)$.
To resolve this issue in Lean, we extend the codomain of the degree function to \texttt{WithBot ($\sigma \to_0 \mathbb{N}$)}. This type adjoins a bottom element $\bot$ that is strictly smaller than every exponent vector. This allows us to separate the degree of the zero polynomial from that of constant polynomials:
\begin{enumerate}
    \item $\deg(0) = \bot,$
    \item $\deg(c) = 0$, where $c$ is  a nonzero constant. 
\end{enumerate}
Under this convention, $\deg(0) < \deg(c)$ is valid and the classical additive properties of degree can be recovered, since $\bot + a = \bot$ for all numbers $a$.
We formalize the degree function with a bottom element as follows.
\begin{leancode}
variable (m) in
def withBotDegree : WithBot (σ →₀ ℕ) :=
    f.support.image m.toSyn |>.max.map m.toSyn.symm
\end{leancode}
Although \texttt{withBotDegree} resolves the ambiguity of the zero polynomial, we sometimes continue to use the existing degree function, since multivariate polynomials and monomial orders are defined over the exponent vector type $\sigma \to_0 \mathbb{N}$ rather than \texttt{WithBot ($\sigma \to_0 \mathbb{N}$)}. Using \lean{degree} therefore allows direct interaction with existing Mathlib APIs without repeatedly eliminating the \texttt{WithBot} constructor or turning $\bot$ to $0$ with \lean{unBot}.
\subsection{Polynomial Division and Remainder}
\label{reduction}
The Gröbner basis theory is fundamentally based on polynomial reduction. We therefore begin by formalizing the notion of remainder.
\begin{definition}[{\ourdef[Remainder]{MonomialOrder.IsRemainder}}]
\label{isremainder}
Let $f$ be a multivariate polynomial and $B$ a set of multivariate polynomials over a field $k$, and fix a monomial order. Suppose that there exists a function $g : B \to \kxi$ with finite support and a polynomial $r$ such that
\begin{itemize}
    \item [(1)] $f = \sum_{b \in B} g(b)\, b + r$,   
    \item [(2)] degree of any $g(b)b$ where $b\in B$ is less than or equal to the degree of $f$.
    \item [(3)] no term of $r$ is divisible by the leading monomial of any nonzero element of $B$. 
\end{itemize}
 Then $r$ is called a remainder of $f$ on division by $B$. We also write $f\xrightarrow{B} r$ for the remainder on division of $f$ by polynomials in  set $B$. 
\end{definition}
\begin{leancode}
def IsRemainder :=
    (∃ (g : B →₀ MvPolynomial σ R),
      f = Finsupp.linearCombination _ (fun (b : B) ↦ b.val) g + r ∧
        ∀ (b : B), m.toWithBotSyn (m.withBotDegree b.val) +
          m.toWithBotSyn (m.withBotDegree (g b)) ≤
          m.toWithBotSyn (m.withBotDegree f)) ∧
            ∀ c ∈ r.support, ∀ b ∈ B, b ≠ 0 → ¬ (m.degree b ≤ c)
\end{leancode} 
Although practical computations typically involve finite families of polynomials, the notion of remainder is most naturally defined with respect to an abstract set of polynomials. To relate these points of view, we establish the following theorem.
\begin{leancode}
theorem isRemainder_range {ι : Type*} (f : MvPolynomial σ R)
    (b : ι → MvPolynomial σ R) (r : MvPolynomial σ R) :
      m.IsRemainder f (Set.range b) r ↔
      (∃ g : ι →₀ MvPolynomial σ R,
        f = Finsupp.linearCombination _ b g + r ∧
        ∀ i : ι, m.withBotDegree (b i * g i) ≼'[m] m.withBotDegree f) ∧
      ∀ c ∈ r.support, ∀ i : ι, b i ≠ 0 → ¬ (m.degree (b i) ≤ c)
\end{leancode}
For practical formalization, we establish a variant of this theorem formulated for finite index sets that simplifies its use in concrete developments.
 
\begin{leancode}
theorem isRemainder_range_fin {ι : Type*} [Fintype ι] (b : ι → MvPolynomial σ R)
    (r : MvPolynomial σ R) :
      m.IsRemainder p (Set.range b) r ↔
      (∃ g : ι → MvPolynomial σ R,
          p = ∑ i : ι, (b i * g i) + r ∧
          ∀ i : ι, m.degree (b i * g i) ≼[m] m.degree p) ∧
        ∀ c ∈ r.support, ∀ i : ι, b i ≠ 0 → ¬ (m.degree (b i) ≤ c)
\end{leancode}
\begin{theorem}[{\ourdef[Existence of Remainders]{MonomialOrder.IsRemainder.exists_isRemainder}}]
\label{remainder_exist}
Fix a monomial order. For any polynomial $f \in \kxi$ and any ordered tuple of polynomials $B = (f_1, \ldots, f_s)$ in $\kxi$, there exist polynomials $a_1, \ldots, a_s \in \kxi$ and a polynomial $r \in \kxi$ such that
\[
f = a_1 f_1 + \cdots + a_s f_s + r,
\]
and $r$ is a remainder of $f$ with respect to the set $\{f_1, \ldots, f_s\}$ in the sense of Definition~\ref{isremainder}. 
\end{theorem}
\begin{leancode}
theorem exists_isRemainder {B : Set (MvPolynomial σ R)}
    (hB : ∀ b ∈ B, IsUnit <| m.leadingCoeff b) (p : MvPolynomial σ R) :
    ∃ (r : MvPolynomial σ R), m.IsRemainder p B r
\end{leancode}
The existence of remainders is a fundamental property underlying polynomial reduction. 
\subsection{Gröbner Basis}\label{sec4_3}
We now turn to the formalization of Gröbner bases in Lean. Building on the notions of multivariate polynomials, monomial orders, and polynomial reduction developed in the previous sections, we introduce the formal definition of Gröbner bases and their fundamental properties.
First, we formalize the notion of leading term defined in Definition~\ref{leadingterm}.
\begin{leancode}
noncomputable def leadingTerm (f : MvPolynomial σ R) : MvPolynomial σ R :=
  monomial (m.degree f) (m.leadingCoeff f)
\end{leancode}
With the notion of leading term in place, we can now define Gröbner bases.
\begin{definition}[{\ourdef[Gröbner Basis]{MonomialOrder.IsGroebnerBasis}}]
\label{gb}
Fix a monomial order  on the polynomial ring $\kxi$. A subset $G$ of an ideal $I \subseteq \kxi$ is said to be a Gröbner basis of $I$ if the leading terms of the elements of $G$ generate the leading term ideal of $I$,
\[
\langle \LT(G) \rangle = \langle \LT(I) \rangle.
\]
where $\LT(G)$ and $\LT(I)$ denote the set of leading terms of elements of $G$ and $I$, respectively.
\end{definition}
  
\begin{leancode}
def IsGroebnerBasis {R : Type*} [CommSemiring R] (G : Set (MvPolynomial σ R))
    (I : Ideal (MvPolynomial σ R)) :=
  G ⊆ I ∧ Ideal.span (m.leadingTerm '' ↑I) = Ideal.span (m.leadingTerm '' G)
\end{leancode}
Having established the definition of Gröbner bases, we turn to the question of their existence. While the definition is formally general, the existence of a finite Gröbner basis for any ideal is fundamentally rooted in the Noetherian property of the polynomial ring. Consequently, the following theorem is stated for polynomial rings in finitely many variables over a field, ensuring that such bases always exist.
\begin{theorem}
\label{exist_gb}
\jumplink{MonomialOrder.IsGroebnerBasis.exists_isGroebnerBasis_finite} Fix a monomial order on $R$, where $R = k[x_1, \dots, x_n]$. Then for any ideal $I \subseteq R$, there exists a finite set $G = \{g_1, \dots, g_t\} \subset I$ such that $G$ is a Gröbner basis of $I$.
\end{theorem}
\begin{leancode}
theorem exists_isGroebnerBasis_finite [Finite σ] :
    ∃ G : Finset (MvPolynomial σ k), IsGroebnerBasis m G ↑I
\end{leancode}
Beyond existence results, Gröbner bases play a central role in solving the ideal membership problem. A key consequence is presented in Theorem~\ref{ideal_mem}.
\begin{theorem}
\label{ideal_mem}
\jumplink{MonomialOrder.IsGroebnerBasis.isRemainder_zero_iff_mem_ideal} Fix a monomial order on the polynomial ring $\kxi$. Let $I$ be an ideal and let $G \subseteq I$. If $G$ is a Gröbner basis of $I$, then for any polynomial $p$, the remainder of $p$ upon division by $G$ is zero if and only if $p \in I$.
\end{theorem}
\begin{leancode}
theorem isRemainder_zero_iff_mem_ideal {p : MvPolynomial σ R}
    {G : Set (MvPolynomial σ R)} {I : Ideal (MvPolynomial σ R)}
    (hG : ∀ g ∈ G, IsUnit (m.leadingCoeff g))
    (h : m.IsGroebnerBasis G I) :
    m.IsRemainder p G 0 ↔ p ∈ I
\end{leancode}
Since every polynomial $p \in I$ reduces to a remainder of zero, it follows that $p$ can be expressed as a linear combination of elements of $G$. This leads to Proposition~\ref{ideal_eq_span}.
\begin{proposition}
\label{ideal_eq_span}
\jumplink{MonomialOrder.IsGroebnerBasis.ideal_eq_span}
Let $I$ be an ideal and $G$ be a subset of the polynomial ring $\kxi$. If $G$ is a Gröbner basis for $I$, then we have  $I=\langle G \rangle.$
\end{proposition}
\begin{leancode}
theorem ideal_eq_span {G : Set (MvPolynomial σ R)}
    {I : Ideal (MvPolynomial σ R)}
    (hG : ∀ g ∈ G, IsUnit (m.leadingCoeff g)) (h : m.IsGroebnerBasis G I) :
    I = Ideal.span G
\end{leancode}
To relate the algebraic definition to the reduction behavior, we establish the following equivalence, characterizing Gröbner bases in terms of ideal generation and reduction. This is a key step toward the remainder-based criterion.
\begin{proposition}
\label{gb_characterization_via_degree}
\jumplink{MonomialOrder.IsGroebnerBasis.isGroebnerBasis_iff_subset_and_degree_le_eq_and_degree_le}
Fix a monomial order on the polynomial ring $\kxi$. Let $I$ be an ideal and let $G \subseteq I$.
Then $G$ is a Gröbner basis of $I$ if and only if $I = \langle G \rangle$ and, for every nonzero
polynomial $p \in I$, there exists $g \in G$ such that $\LM(g) \mid \LM(p)$.
\end{proposition}
\begin{leancode}
theorem isGroebnerBasis_iff_subset_and_degree_le_eq_and_degree_le (G : Set (MvPolynomial σ R))
    (I : Ideal (MvPolynomial σ R)) (hG : ∀ g ∈ G, IsUnit (m.leadingCoeff g)) :
    m.IsGroebnerBasis G I ↔
      Ideal.span G = I ∧ ∀ p ∈ I, p ≠ 0 → ∃ g ∈ G, m.degree g ≤ m.degree p
\end{leancode}
 
The following proposition characterizes the Gröbner bases in terms of remainders and provides a criterion for deciding ideal membership.
\begin{proposition}
\label{gb_remainder}
\jumplink{MonomialOrder.IsGroebnerBasis.isGroebnerBasis_iff_subset_ideal_and_isRemainder_zero}
Fix a monomial order on the polynomial ring $\kxi$. Let $I$ be an ideal of $\kxi$. A set $G$ is a Gröbner basis of $I$ if and only if $G \subseteq I$ and, for every $f \in I$, the remainder of $f$ upon division by $G$ is $0$.
\end{proposition}
\begin{leancode}
theorem isGroebnerBasis_iff_subset_ideal_and_isRemainder_zero
    (G : Set (MvPolynomial σ R)) (I : Ideal (MvPolynomial σ R))
    (hG : ∀ g ∈ G, IsUnit (m.leadingCoeff g)) :
    m.IsGroebnerBasis G I ↔ G ⊆ I ∧ ∀ p ∈ I, m.IsRemainder p G 0
\end{leancode}
 
Finally, an important property of Gröbner bases is that they provide a unique remainder for elements modulo the ideal they generate. Unlike a general multivariate division, where the remainder typically depends on the order of the divisors, division by a Gröbner basis guarantees  a unique remainder.
\begin{theorem}
\label{unique_of_remainder}
\jumplink{MonomialOrder.IsGroebnerBasis.existsUnique_isRemainder}
Fix a monomial order on the polynomial ring $\kxi$. Let $I$ be an ideal and let $G = \{g_1, \dots, g_t\}$ be a Gröbner basis of $I$. Then for any polynomial $p$, the remainder of $p$ upon division by $G$ is unique.
\end{theorem}
\begin{leancode}
theorem existsUnique_isRemainder {G : Set (MvPolynomial σ R)}
    {I : Ideal (MvPolynomial σ R)}
    (h : m.IsGroebnerBasis G I) 
    (hG : ∀ g ∈ G, IsUnit (m.leadingCoeff g)) (p : MvPolynomial σ R) :
    ∃! (r : MvPolynomial σ R), m.IsRemainder p G r
\end{leancode}
\subsection{Buchberger Criteria}\label{sec4_buchberger}
To formally verify Buchberger's criterion in Lean~4, we first introduce $S$-polynomials.
Our definition of the $S$-polynomial differs slightly from the standard formulations in~\cite{becker1993grobner,cox2025ideals}. We adopt this variant to obtain a more natural implementation in Lean, in particular avoiding the need to assume invertibility of leading coefficients.
\begin{definition}[{\ourdef[$S$-polynomial]{MonomialOrder.sPolynomial}}]\label{spoly}
Fix a monomial order. Let $f,g \in \kxi$. Define the $S$-polynomial by
\[
  \spoly(f, g) = \frac{\LC(g)x^\gamma}{\mathrm{LM}(f)} \cdot f - \frac{\LC(f)x^\gamma}{\mathrm{LM}(g)} \cdot g,
\]
 where $x^\gamma$ denotes the least common multiple of $\LM(f)$ and $\LM(g)$.
\end{definition}
Since this formulation differs from the standard definition only by multiplication by a nonzero scalar factor arising from clearing leading coefficients,
\[
\spoly(f,g)=\frac{x^\gamma}{\LT(f)}\,f-\frac{x^\gamma}{\LT(g)}\,g,
\]
it does not affect the validity of Buchberger's criterion.
We now give the formal definition of the $S$-polynomial following Definition~\ref{spoly}.
\begin{leancode}
noncomputable def sPolynomial (f g : MvPolynomial σ R) : MvPolynomial σ R :=
    monomial (m.degree g - m.degree f) (m.leadingCoeff g) * f -
    monomial (m.degree f - m.degree g) (m.leadingCoeff f) * g
\end{leancode}
Before proving Buchberger's criterion, we first establish a lemma characterizing the cancelation of leading terms. If two nonzero polynomials have the same leading degree, then the leading degree of their $S$-polynomial is strictly smaller. More generally, if a sum of polynomials with the same leading degree exhibits cancelation of leading terms, then the sum can be expressed as a linear combination of their $S$-polynomials with coefficients in $k$.

\begin{leancode}
lemma sPolynomial_decomposition_of_degree_sum_smul_le₀ {R} [CommRing R] {d : m.syn} {ι : Type*}
    {B : Finset ι} {c : ι → R} {g : ι → MvPolynomial σ R}
    (hd : ∀ b ∈ B,
      (m.toSyn <| m.degree <| g b) = d ∧ IsUnit (m.leadingCoeff <| g b) ∨ g b = 0)
    (hfd : (m.toSyn <| m.degree <| ∑ b ∈ B, c b • g b) < d) :
    ∃ (c' : ι → ι → R),
      ∑ b ∈ B, c b • g b = ∑ b₁ ∈ B, ∑ b₂ ∈ B, (c' b₁ b₂) • m.sPolynomial (g b₁) (g b₂)
\end{leancode}
Then we attempt to formalize Buchberger's criterion, which characterizes Gröbner bases via the reduction of $S$-polynomials of critical pairs. This replaces a universal quantifier on an infinite set $\langle G\rangle$ by finitely many checks relative to $G$.
\begin{theorem}[{\ourdef[Buchberger's Criterion]{MonomialOrder.IsGroebnerBasis.isGroebnerBasis_iff_isRemainder_sPolynomial_zero}}]
\label{buchberger}
Fix a monomial order on the polynomial ring $\kxi$. A set $G$ is a Gröbner basis of the ideal $\langle G\rangle$ if and only if for all $f,g \in G$, $0$ is the remainder of $S$-polynomial $\spoly(f,g)$ on division of $G$.
\end{theorem}
\begin{runningexample}\label{eg_buchberger}
Consider the Cyclic--3 ideal $I \subset \mathbb{Q}[x_0,x_1,x_2]$ with lexicographic order $x_0 > x_1 > x_2$, generated by
\[
\{\,f_1=x_0+x_1+x_2,\quad f_2=x_0x_1+x_1x_2+x_2x_0,\quad f_3=x_0x_1x_2-1\,\}.
\]
A Gröbner basis for $I$ is 
\begin{equation}\label{reducedG}
G=\{\,g_1=x_0+x_1+x_2,\quad g_2=x_1^2+x_1x_2+x_2^2,\quad g_3=x_2^3-1\,\}.
\end{equation}
By Buchberger's criterion, it suffices to check that for all pairs $g_i,g_j\in G$, the $S$-polynomial $\spoly(g_i,g_j)$ reduces to $0$ modulo $G$. As an illustration, consider $(g_1,g_2)$. We have $\LM(g_1)=x_0$, $\LM(g_2)=x_1^2$, and
$\mathrm{lcm}(\LM(g_1),\LM(g_2))=x_0x_1^2$, hence
\begin{align*}
\spoly(g_1,g_2)
&= \frac{x_0x_1^2}{x_0}\,g_1 - \frac{x_0x_1^2}{x_1^2}\,g_2 \\
&= x_1^2(x_0+x_1+x_2) - x_0(x_1^2+x_1x_2+x_2^2) \\
&= -x_0x_1x_2 - x_0x_2^2 + x_1^3 + x_1^2x_2\\
&=(-x_1x_2-x_2^2)g_1 + (x_1+x_2)g_2.
\end{align*}
so $\spoly(g_1,g_2)$ reduces to $0$ modulo $G$. The remaining pairs are checked similarly.
\end{runningexample}
We now formalize Buchberger's criterion, as stated in Theorem~\ref{buchberger}. 
\begin{leancode}
theorem isGroebnerBasis_iff_ideal_eq_span_and_isGroebnerBasis_span₀ {R} 
    [CommRing R] {I : Ideal (MvPolynomial σ R)} {G : Set (MvPolynomial σ R)}
    (hG : ∀ g ∈ G, IsUnit (m.leadingCoeff g) ∨ g = 0) :
    m.IsGroebnerBasis G I ↔ I = Ideal.span G ∧
      ∀ (g₁ g₂ : G), m.IsRemainder (m.sPolynomial g₁ g₂ : MvPolynomial σ R) G 0
\end{leancode}
\subsection{Reduced Gröbner Basis}\label{sec4_reduced_gb}
Since the definition of Gröbner basis only requires that the leading terms of its elements generate the leading term ideal $\langle \LT(I) \rangle$, redundant generators may occur. For example, if $g_1, g_2 \in G$,   $\LM(g_1)$ divides $\LM(g_2)$, then $g_2$ does not contribute any new information. The notion of a minimal Gröbner basis is introduced to eliminate such redundant generators.
\begin{definition}[{\ourdef[Minimal Gröbner Basis]{MonomialOrder.IsGroebnerBasis.IsMinimal}}]
\label{minimal_gb}
Fix a monomial order on the polynomial ring $\kxi$. A subset $G \subseteq I$ is called a minimal Gröbner basis of the ideal $I$ if the following conditions hold:
\begin{itemize}
    \item[(1)] $G$ is a Gröbner basis of $I$.
    \item[(2)] Every $p \in G$ is monic.
    \item[(3)] For all distinct $p,q \in G$, $\LM(p) \nmid \LM(q)$.
\end{itemize}
\end{definition}
We now give the Lean formalization of the notion of a minimal Gröbner basis, following Definition~\ref{minimal_gb}:
\begin{leancode}
def IsGroebnerBasis.IsMinimal (hG : m.IsGroebnerBasis G I) :=
  (∀ p ∈ G, m.Monic p) ∧
  (∀ p ∈ G, ∀ q ∈ G, q ≠ p → ¬ m.degree q ≤ m.degree p)
\end{leancode}
\begin{runningexample}[continued]\label{eg_minimal_gb}
One can verify that  $G=\{g_1,g_2,g_3\}$ (\ref{reducedG}) satisfies the conditions of Definition \ref{minimal_gb}, and therefore is a minimal Gröbner basis of $I$.
\end{runningexample}
While a minimal Gröbner basis removes redundant generators, it is not necessarily unique for a fixed ideal and monomial order. Different minimal Gröbner bases may exist for the same ideal. To obtain a canonical representation, we impose a stronger condition: no monomial appearing in any polynomial of the basis is divisible by the leading monomial of any other polynomial in the set. This leads to the notion of a reduced Gröbner basis.
\begin{definition}[{\ourdef[Reduced Gröbner Basis]{MonomialOrder.IsGroebnerBasis.IsReduced}}]
\label{reduced_gb}
Fix a monomial order on the polynomial ring $\kxi$. A subset $G \subseteq I$ is called a reduced Gröbner basis of the ideal $I$ if the following conditions hold:
\begin{itemize}
    \item[(1)] $G$ is a Gröbner basis of $I$.
    \item[(2)] Every $p \in G$ is monic.
    \item[(3)] For all distinct $p,q \in G$, no monomial appearing in $p$ is divisible by $\LM(q)$.
\end{itemize}
\end{definition}
We now formalize the notion of a reduced Gröbner basis in Lean, following Definition~\ref{reduced_gb}:
\begin{leancode}
def IsGroebnerBasis.IsReduced (hG : m.IsGroebnerBasis G I) :=
  (∀ p ∈ G, m.Monic p) ∧ ∀ p ∈ G, m.IsRemainder p (G \ {p}) p
\end{leancode}
\begin{runningexample}[continued]\
It is straightforward to check that
\[
G'=\{g'_1,g_2,g_3\}, \qquad 
g'_1 = x_0 + x_1^2 + x_1 x_2 + x_1 + x_2^2 + x_2,
\]
is also a minimal Gröbner basis of the same ideal \(I\).
However, \(G'\) is not reduced, since \(g'_1\) contains the monomial \(x_1^2\), which is divisible by \(\LM(g_2)=x_1^2\).
Reducing \(g'_1\) by \(g_2\), namely setting \(g_1 = g'_1 - g_2\), removes this reducible term and yields the reduced Gröbner basis \(G=\{g_1,g_2,g_3\}\) (see~\eqref{reducedG}).
    
\end{runningexample}
The stronger condition in Definition~\ref{reduced_gb} yields the following fundamental property of reduced Gröbner bases.
\begin{theorem}
\label{uniqueExists_reduced_gb}
\jumplink{MonomialOrder.IsGroebnerBasis.IsReduced.uniqueExists}
Fix a monomial order on the polynomial ring $\kxi$. Every ideal $I \subseteq \kxi$ admits a unique reduced Gröbner basis.
\end{theorem}
\begin{leancode}
theorem IsReduced.uniqueExists {k} [Field k] (I : Ideal (MvPolynomial σ k)) :
    ∃! (B : Set (MvPolynomial σ k)), ∃ (h : m.IsGroebnerBasis B I), h.IsReduced
\end{leancode}
The following version is in the general commutative ring instead of field, with the extra hypothesis that the leading coefficients of a Gröbner basis are units.
\begin{leancode}
theorem IsReduced.uniqueExists_of_isGroebnerBasis {R} [Nontrivial R] [CommRing R]
    {G : Set (MvPolynomial σ R)} {I : Ideal (MvPolynomial σ R)}
    (hG : m.IsGroebnerBasis G I) (hG' : ∀ g ∈ G, IsUnit (m.leadingCoeff g)) :
    ∃! (B : Set (MvPolynomial σ R)), ∃ (h : m.IsGroebnerBasis B I),
      h.IsReduced
\end{leancode}
\section{Finite Characterization of Infinite Gröbner Bases}
\label{infinite_to_finite}
In this section, we explain how Gröbner bases in polynomial rings with infinitely many variables can be characterized via finite-variable reductions. Let $\sigma$ and $\sigma'$ be variable index types, and write two polynomial rings 
\[
S = \kxi, \qquad S' = \kxiv.
\]
Let $\M := (\sigma \to_0 \mathbb{N})$ and $\M' := (\sigma' \to_0 \mathbb{N})$ be the corresponding types of exponent vectors, and denote monomial orders on $\M$ and $\M'$ as $(\M, \le_\M)$ and $(\M', \le_{\M'})$, respectively.
A key prerequisite for relating these two settings is the compatibility of monomial orders. We now introduce the notion of an embedding of monomial orders.
\begin{definition}[{\ourdef[Monomial Order Embedding]{MonomialOrder.Embedding}}]
\label{order_embedding}
Given two monomial orders $(M, \le_{m})$ and $(M', \le_{m'})$, a function $f : \sigma' \to \sigma$ is an order embedding if $f$ satisfies 
\begin{enumerate}
    \item[(1)] The function $f$ is injective.
    \item[(2)] Denote the lift from $\M'$ to $\M$ via $f$ as $f' : \M' \to \M$.
    For all $u$, $v$ in $\M'$, if $u\le_{m'}v$, then $f'(u)\le_mf'(v)$.
\end{enumerate}
\end{definition}
We formalize Definition~\ref{order_embedding} in Lean as follows:
\begin{leancode}
structure Embedding extends σ' ↪ σ where
  monotone' : Monotone (m.toSyn ∘ Finsupp.mapDomain toFun ∘ m'.toSyn.symm)
\end{leancode}
Given a monomial order on $\M$, any injective map $f:\sigma' \to \sigma$ induces a monomial order on $\M'$. This order is obtained by pulling back the order on $\M$ along the embedding of exponent vectors induced by $f$
\begin{leancode}
open ofInjective in
noncomputable def ofInjective {σ' : Type*} {f : σ' → σ} (hf : f.Injective) :
    MonomialOrder σ' :=
  { syn := Syn m f
    toSyn : (σ' →₀ ℕ) ≃+ (Syn m f) := toSyn' m f
    ... }
\end{leancode}
With this construction in place, the induced monomial order on $\M'$ admits a natural embedding into the monomial order on $\M$.
\begin{leancode}
def Embedding.ofInjective {f : σ' → σ} (hf : f.Injective) :
    Embedding (m.ofInjective hf) m where
  toEmbedding := ⟨f, hf⟩
  monotone' := ...
\end{leancode}
\begin{lemma}
\jumplink{MonomialOrder.Embedding.isRemainder_killCompl_of_isRemainder_rename}
Let $p \in S'$, $r \in S$, and let $B \subseteq S'$ be a set of polynomials. If $r$ is a remainder of $p$ on division by $B$ in $S$, then the restriction $r'$ of $r$, obtained by setting variables outside $\sigma'$ to $0$, is a remainder of $p$ on division by $B$ in $S'$.
\end{lemma}
\begin{leancode}
lemma MonomialOrder.Embedding.isRemainder_killCompl_of_isRemainder_rename
    {σ' σ} {m' : MonomialOrder σ'} {m : MonomialOrder σ}
    (e : Embedding m' m) {p : MvPolynomial σ' R} {B : Set (MvPolynomial σ' R)}
    {r : MvPolynomial σ R} (h : m.IsRemainder (p.rename e) (rename e '' B) r) :
    m'.IsRemainder p B (r.killCompl e.coe_injective)
\end{leancode}
The following lemma shows that the remainders are invariant under embedding.
\begin{lemma}
\label{remainder_invar}
\jumplink{MonomialOrder.Embedding.isRemainder_iff_isRemainder_rename}
Let $p, r \in S'$ and let $B \subseteq S'$ be a set of polynomials. Then $r$ is a remainder of $p$ on division by $B$ over $S'$ if and only if $r$ is a remainder of $p$ on division by $B$ over $S$.
\end{lemma}
\begin{leancode}
lemma MonomialOrder.Embedding.isRemainder_iff_isRemainder_rename {σ' σ}
    {m' : MonomialOrder σ'} {m : MonomialOrder σ}
    (e : Embedding m' m) (p : MvPolynomial σ' R) (B : Set (MvPolynomial σ' R))
    (r : MvPolynomial σ' R) :
    m'.IsRemainder p B r ↔ 
      m.IsRemainder (p.rename e) (MvPolynomial.rename e '' B) (r.rename e) 
\end{leancode}
There is also invariance on Gr\"{o}bner bases.
\begin{lemma}
\label{gb_invar}
\jumplink{MonomialOrder.Embedding.isGroebnerBasis_iff_isGroebnerBasis_rename}
Let $G \subseteq S'$ be a set of polynomials and let $I \subseteq S'$ be an ideal. Then $G$ is a Gröbner basis for $I$ over $S'$ if and only if $G$ is a Gröbner basis for the ideal generated by $I$ over $S$.
\end{lemma}
\begin{leancode}
lemma MonomialOrder.Embedding.isGroebnerBasis_iff_isGroebnerBasis_rename {σ'}
    {m' : MonomialOrder σ'} {m : MonomialOrder σ}
    (e : Embedding m' m) (G : Set (MvPolynomial σ' R))
    (I : Ideal (MvPolynomial σ' R)) :
    m'.IsGroebnerBasis G I ↔
      m.IsGroebnerBasis (MvPolynomial.rename e '' G) (I.map (MvPolynomial.rename e)) 
\end{leancode}
In general, a finite Gröbner basis need not exist for an ideal in a polynomial ring with infinitely many variables, since such rings are not Noetherian. However, using the lemmas above, one can obtain a finite Gröbner basis by restricting to the subring generated by the variables appearing in a finite generating set, which is Noetherian. This basis can then be lifted to the original polynomial ring. In particular, an ideal is finitely generated if and only if it admits a finite Gröbner basis.
\begin{lemma}
\label{gb_exist'}
\jumplink{MonomialOrder.IsGroebnerBasis.exists_isGroebnerBasis_finite_of_exists_span_finite}
Let $I \subseteq \kxi$ be a finitely generated ideal. Then $I$ admits  a finite Gröbner basis.
\end{lemma}
\begin{leancode}
lemma exists_isGroebnerBasis_finite_of_exists_span_finite {B : Set (MvPolynomial σ k)} 
(hB : B.Finite) :
    ∃ G : Set (MvPolynomial σ k), m.IsGroebnerBasis G (Ideal.span B) ∧ G.Finite
\end{leancode}
In \cite[Theorem 1.12]{IimaYoshino08}, Iima and Yoshino characterize reduced Gröbner bases over polynomial rings with countably many variables via the limit inferior of reduced Gröbner bases on finite-variable subrings. We generalize this result to arbitrary index types using limit inferior constructions over filters.
\begin{theorem}
\label{reduced_gb_lim}
\jumplink{MonomialOrder.IsGroebnerBasis.IsReduced.isReduced_liminf}
Fix an ideal $I \subseteq \kxi$, an index set $\iota$, let  $F \subseteq \mathcal{P}(\iota)$ be a filter over lattice $(\mathcal{P}(\iota), \subseteq)$ s.t. $\emptyset \notin F$. Assume for all $j \in \iota$, $\sigma_{j} \subseteq \sigma$ and $G_j$ be the reduced Gr\"{o}bner basis of ideal $I \cap k[x_i]_{i\in\sigma_j}$ on $k[x_i]_{i\in\sigma_j} \subseteq \kxi$. If $\bigcup_{s \in F}\bigcap_{j \in s} \sigma_j = \sigma$ , then
\[ \bigcup_{s \in F}\bigcap_{j \in s} G_j \]
is the reduced Gr\"{o}bner basis of $I$.
\end{theorem}
\begin{leancode}
theorem IsReduced.isReduced_liminf {k} [Field k] {I : Ideal (MvPolynomial σ k)} 
    {α} {σ' : α → Type*} {m' : (a : α) → MonomialOrder (σ' a)}
    {f : Filter α} [inst_neBot : f.NeBot]
    {e : (a : α) → (m' a).Embedding m}
    (hI : Set.univ = f.liminf (fun x ↦ Set.range (e x)))
    {G' : (a : α) → Set (MvPolynomial (σ' a) k)}
    (hG' : ∀ a, (m' a).IsGroebnerBasis (G' a) (I.comap <| rename (e a)))
    (hG'' : ∀ a, (hG' a).IsReduced) :
    ∃ h : m.IsGroebnerBasis (f.liminf fun a ↦ rename (e a) '' G' a) I, h.IsReduced
\end{leancode}
In particular, letting $\sigma$ and $\iota$ be $\mathbb{N}$, $F$ be the cofinite filter on $\mathbb{N}$, $\sigma_n$ be $[0, n] \cap \mathbb{N}$, and $G_n$ be the reduced Gr\"{o}bner basis of $I\cap k[x_i]_{i\in [0, n] \cap \mathbb{N}}$, then 
\[\bigcup_{\substack{n \in \mathbb{N},\\ [n, +\infty) \subseteq A\subseteq \mathbb{N}}}\bigcap_{i \in A} G_i = \bigcup_{n \in \mathbb{N}}\bigcap_{i\ge n}G_i
\] is the reduced Gr\"{o}bner basis of $I$ \cite[Theorem 1.12]{IimaYoshino08}.
\section{Conclusion and Future Work}\label{future}
Currently, our formalization is based on \lean{MvPolynomial} from Mathlib. This structure provides strong mathematical generality by supporting polynomial rings indexed by arbitrary types, but is not designed for efficient executable algebraic computation. Consequently, our present development focuses on theoretical verification, and in this work, we do not formalize an executable implementation of Buchberger's algorithm.
Future work will address the gap between abstract theory and concrete computation. In particular, we plan to support certificate-based verification of polynomial identities by combining external computer algebra systems with our formal theory. In this workflow, Gröbner bases are computed externally and then verified inside Lean using our formalization of Buchberger's criterion. This approach avoids the computational limitations inherent in working directly with \lean{MvPolynomial}, while preserving full machine-checking.
\section*{Acknowledgment}\label{sec6}
Hao Shen, Junqi Liu, and Lihong Zhi are supported by the National Key R\&D Program of China (Grant No.~2023YFA1009401). Junyu Guo is encouraged and supported by his supervisor, Xishun Zhao, to work on this project.
We sincerely thank several members of the Lean community who reviewed our formalization submitted to Mathlib and provided valuable feedback and concrete improvements. In alphabetical order, we thank Dagur Asgeirsson, Riccardo Brasca, Snir Broshi, Thomas Browning, Bryan Gin-ge Chen, Johan Commelin, Jovan Gerbscheid, Aaron Liu, Jireh Loreaux, Bhavik Mehta, Kyle Miller, Ruben Van de Velde, Eric Wieser, and Andrew Yang for their generous help and support. In particular, Snir Broshi suggested formalizing Gröbner bases using \mathlibdef{Set} rather than \mathlibdef{Finset}, which led us to consider properties specific to infinite Gröbner bases. There are more contributors who have built the infrastructure of Lean ecosystem we rely on, including but not limited to Antoine Chambert-Loir who formalized the properties of the divsion of multivariate polynomials, which we use directly to prove theorem \ref{remainder_exist}.
 
\bibliography{ARL}
\end{document}